\magnification=1200
\font\rmtwelve=cmbx10 at 12pt
\hsize=6truein 
\voffset=.9truein
\vsize=8 true in
\baselineskip=12pt plus 2pt minus 2pt
\parskip=3pt
\tolerance= 10000
\centerline{\rmtwelve  Symplectic Parshin-Arakelov inequality}
\bigbreak
\bigskip
\centerline{Tian-Jun Li}
\bigskip
\noindent{\bf  \S1. Introduction}
\medskip

Lefschetz fibration is the smooth analogue of stable holomorphic 
fibration. In dimension four, its importance stems from its close 
relations to the mapping class groups and the Deligne-Mumford 
moduli space of algebraic curves. Recently it has received wide 
attention because of the discovery, mainly due to Donaldson, that 
it provides a topological description of symplectic manifolds.

For a relatively minimal genus $g$ stable holomorphic fibrations
over a genus $h>0$  Riemann surface, there is a famous
Parshin-Arakelov inequality: $c_1^2\geq 8(g-1)(h-1)$. 
In this paper, we will present its symplectic analogue.

\noindent{\bf Theorem 1}. Let $M$ be a relatively minimal genus 
$g$ Lefschetz fibration over a genus $h$ surface. If $M$ is not
rational or ruled, then $$c_1^2(M)\geq 2(g-1)(h-1).$$ and it is 
sharp in the case $h=0$.

When $h$ is positive, the inequality generalizes Kotschick's result 
[K] for surface bundles. We do not know whether it is sharp or not.

In the theorem, the condition  
that $M$ not being rational or ruled is necessary. Because
when $h=0$, our inequality is $c_1^2\geq 2-2g$,
while there are many Lefschetz fibrations over $S^2$ on rational and 
ruled surfaces with $c_1^2= 4-4g$.  Rational and ruled 
surfaces are symplectic four-manifolds with exceptional properties 
and can be characterized among all symplectic four manifolds 
in several ways (see [L], [Liu], [Mc]).  In this paper, we also 
analyze Lefschetz fibrations on these manifolds.  The analysis 
of Lefschetz fibrations on ruled surfaces, in conjuction with 
Theorem 1, allows us to obtain a lower bound of the number of 
irreducible singular fibers for Lefschetz fibrations over $S^2$.

\noindent{\bf Theorem 2}. The number of irreducible 
singular fibers of a genus $g$ Lefschetz fibration over $S^2$
is no less than $g$.

%\noindent 1. $n\geq g$ if $g$ is even and $n\geq g+1$ if $g$ is odd.

%\noindent 2. $n\geq 4$ is $g=2$ and $n\geq 6$ if $g=3$.

The organization of this paper is as follows. We review Lefschetz
fibrations in \S2. In \S3, we present the proof of the symplectic
Parshin-Arakelov inequality. In \S4, we first study Lefschetz 
fibrations on ruled surfaces. We then present the estimate of the 
minimal number of irreducible singular fibers. Finally we also discuss
the applications of our theorems to the mapping class groups and the Deligne-Mumford 
moduli space of algebraic curves.

The author thanks A. Stipsicz for many stimulating discussions. 
He also wishes to thank K. F. Liu, F. Luo, Y. Matsumoto, B. Ozbagci, 
T. Pantev, I. Smith, G. Tian and S. T. Yau for their interest in this work.
He is grateful to the referee, whose many  suggestions greatly improved
the quality of the presentation.
This research is partially supported by NSF DM-9975469.

\bigskip
\noindent{\bf \S2.  Lefschetz fibrations}
\medskip

\noindent{\bf Definition 2.1}. Let $M$ be a compact, connected,
oriented smooth four-manifold. A Lefschetz fibration is a map
$\pi:M\longrightarrow \Sigma$, where $\Sigma$ is a compact, connected,
oriented surface and $\pi^{-1}(\partial \Sigma)=\partial M$,
such that 

\noindent a). the set of critical points $C=\{x_1, \cdots, x_n\}$
of $\pi$ is non-empty and lies in interior of $M$; 

\noindent b). about each $x_i$ and $\pi(x_i)$, there are 
orientation-preserving complex local coordinate charts on which
$$\pi (z_1, z_2)=z_1^2+z_2^2;$$

\noindent c). $\pi$ is injective on $C$.

A regular fiber is a closed smooth surface, its genus called the genus 
of the Lefschetz fibration. Each singular fiber is a transversely 
immersed surface with a positive double point. A singular fiber is 
called reducible if the connected component containing the critical 
point becomes disconnected after the critical point is removed. 
b) and c) imply that a reducible fiber has exactly two components,
each with square $-1$. 
A Lefschetz fibration is relatively minimal if there is no singular
fiber containing a sphere of self-intersection $-1$.

The existence of a Lefschetz fibration $\pi:M\longrightarrow \Sigma$ 
with regular fiber $F$ provides a handlebody description of $M$ 
(see [K] for more details). A regular neighborhood of a singular fiber
is diffeomorphic to $(F\times D^2)\cup H^2$, where $H^2$ is a 
$2-$handle attached along a simple closed curve $\gamma$ in a fiber 
$F\times \{pt\}$ in the boundary of $F\times D^2$. The attaching 
circle $\gamma$, well-defined up to isotopy, is called the vanishing 
cycle. The boundary of $(F\times D^2)\cup H^2$ is diffeomorphic 
to a $F-$bundle over $S^1$ whose monodromy is given by the 
right-handed Dehn twist about $\gamma$, $D(\gamma)\in {\cal MC}(F)$,
where ${\cal MC}(F)$ is the mapping class group of $F$. Geometrically,
as one approaches the singular fiber, the vanishing cycle is
shrunk to the critical point. We see that a separating vanishing
cycle corresponds to a reducible fiber.

If $\Sigma$ is a two disc, then $M$ is diffeomorphic to 
$(F\times D^2)\cup H^2_i\cup\cdots\cup H^2_n$, where each 
two-handle $H^2_i$ is attached along a vanishing cycle 
$\gamma_i$ in a fiber in the boundary of $\Sigma\times D^2$. 
The boundary of this Lefschetz fibration is a $F-$bundle 
over $S^1$ whose monodromy is the product $D(\gamma_1)\cdots 
D(\gamma_n)$. When $\Sigma$ is $S^2$, we get a Lefschetz 
fibration $M_0$ over $D^2$ by removing a regular neighborhood 
$U$ from a regular fiber. Since $U$ is a trivial $F-$bundle 
over a two disc, the boundary of $M_0$ must be a trivial 
$F-$bundle over $S^1$. The global monodromy $D(\gamma_1)\cdots 
D(\gamma_n)$ is therefore trivial.  $M$ can thus be described as 
$$M=(F\times D^2)\cup H^2_1\cup\cdots\cup H^2_n \cup (F\times D^2). 
\eqno(1)$$ The converse is also true: a relator $D(\gamma_1)\cdots
D(\gamma_n)=1$ gives rise to a Lefschez fibration over $S^2$.

It is not difficult to prove that Lefschetz fibration of genus zero 
must be a blow-up of $S^2-$bundle over a closed surface. Genus one 
Lefschetz fibrations are also well understood thanks to the work
of Kas, Moishezon, Mandelbaum, Harper and Matsumoto (see [M1]). The
relatively minimal ones are fiber sums of torus bundles and $E(1)$.

Recently, Donaldson [D] obtained a remarkable result concerning
the existence of Lefschetz fibrations. Before stating Donaldson's 
result, let us first introduce the definition of a symplectic 
Lefschetz fibration.

\noindent{\bf Definition 2.2}. A Lefschetz fibration
$M\longrightarrow \Sigma$ is called a symplectic 
Lefschetz fibration if there exists a symplectic form $\omega$
on $M$, such that for any $p\in \Sigma$, $\omega$ is nondegenerate
at each smooth point on the fiber $F_p$, and that at each double point,
$\omega$ is nondegenerate on the two planes contained in the
tangent cone. 

Donaldson proves that any symplectic
four-manifold admits symplectic Lefschetz fibrations over $S^2$
after perhaps blowing up. Gompf proves that (see also [ABKP], 
[ST]) most Lefschetz fibrations admit a symplectic structure.

\noindent{\bf Thereom 2.3} ([GS]). If a four-manifold admits 
a Lefschetz fibration $\pi:M\longrightarrow \Sigma$ and the 
fiber represents an essential class, then $M$ admits a 
symplectic Lefschetz fibration structure. In particular, 
when $g\geq 2$, $M$ admits a symplectic Lefschetz fibration 
structure.

We now give some elementary lemmas for Lefschetz fibrations over $S^2$, 
which will be used in \S4.

\noindent{\bf Lemma 2.4}. Let $M\longrightarrow S^2$ be a Lefschetz 
fibration over $S^2$ with regular fiber $F$. 
Let $l$,  $s$ and $n$ be the number 
of singular fibers, reducible singular fibers and irreducible 
singular fibers
respectively. Then,
 
\noindent 1. $n\geq b_1(F)-b_1(M)$, and $n=0$ iff $b_1(F)=b_1(M)$;

\noindent 2. $s+1\leq b^-\leq l+1$, $1\leq b^+\leq n+1$;

\noindent 3. $\sigma=4k-l$ for some non-negative integer $k$;
if all the singular fibers are reducible, then $\sigma=-l$.

\noindent{\it Proof}. Part 1 is well known since 
non-separating vanishing cycles represent nontrivial classes in
$H_1(F)$ and, from
the handlebody description, they 
generate the kernel of the natural map $H_1(F)\longrightarrow
H_1(M)$ induced by inclusion.  We first prove part 2.
By Gompf's theorem, $M$ has symplectic structure, and so $b^+\geq 1$. 
Since $b_2$ is bounded by $l+2$, we immediately get the upper 
bound for $b^-$. To show $b^-\geq s$, let $G_1, \cdots, G_s$ be the 
connected components of each reducible singular fiber. We know that
$G_1^2=\cdots=G_s^2=-1$ and $G_i\cdot G_j=0$ for $i\ne j$.
Thus the intersection form on the subspace generated by
$G_1,\cdots, G_s$ is negative definite. A regular fiber 
$F$ is orthogonal to this $s-$dimensional subspace and
has square zero.  Thus $b^-\geq s+1$ because the intersection 
form is nondegenerate. Since $b_2\leq l+2$, the upper bound of $b^+$
follows.  

Now we turn to  part 3. By the handlebody description of a Lefschetz 
fibration, there are $l+2$ two-handles, $2g$ one-handles and 
$2g$ three-handles. By Poincare Duality, $l+2-b_2=2(2g-b_1)$.
Since $b_2=2b^+-\sigma$ and $b^+\equiv b_1-1\pmod 2$, we find that
$l\equiv -\sigma\pmod 4.$  If all the singular fibers are 
reducible, i.e. $s=l$, then $b^-$ must be $l+1$ and $b^+=1$.  
Therefore $\sigma=-l$.

Lemma 2.4 can also be proved using the signature computation in [O].

\noindent{\bf Lemma 2.5}. For any genus $g$ Lefschetz fibrations with
$\sigma\geq -l+4$, $b_1, b_2, b^+$ and $\sigma$ have upper bounds
$ 2g-2, l-2, n-3$ and $ n-s-4$ respectively.

\noindent{\it Proof}. Since we assume that $\sigma\geq -l+4$, 
there exists irreducible singular fibers by  part 3
of Lemma 2.4, i.e. $n>0$.
So $b_1=2g$ is impossible because it would imply that $n=0$ 
by  part 1 of Lemma 2.4.
%hence all the singular fibers are reducible. But such Lefschetz fibrations do %not exist by the third part of Lemma 2.4.
If $b_1=2g-1$, by  part 1 of 
Lemma 2.4,  the non-separating vanishing cycles 
generate a rank one subgroup of $H_1(F)$. But this is again
impossible since the action of a 
Dehn twist along a non-separating curve on $H_1(F)$ is of 
infinite order.  Thus we have shown that $b_1\leq 2g-2$.  

This upper bound of $b_1$, plus  the handlebody description,
implies $b_2\leq (l+2)- 2\cdot 2=l-2$.  
Finally this upper bound of $b_2$ gives
the upper bounds of $b_2^+$ and $\sigma$  with  part 2 of Lemma 2.4.

To end this section, we describe the connection
between Lefschetz fibrations and the Deligne-Mumford
moduli space of stable curves $\overline{\cal M}_g$ (see [Sm2]). 
Recall that for $g\geq 2$,
$\overline{\cal M}_g$ is the stable compactification of 
${\cal M}_g$, the moduli space of curves of genus $g$.
It is a projective orbifold, and the compactifying divisor
${\cal C}=\overline{\cal M}_g-{\cal M}_g$ consists of stable curves
with at least one node.

By choosing a metric compatible with the sympelctic form and K\"ahler
in the neighborhood of the singular fibers, we can
obtain a smooth map from the two-sphere to $\overline{\cal M}_g$.
This map is restricted to intersect with ${\cal C}$, and 
each intersection point is transverse,
positive, and lies outside the locus of curves with more than
one node.
Such a map is well defined up to isotopy preserving the
condition on the intersection with ${\cal C}$.
Stable Kahler fibrations correspond to holomorphic maps.
On $\overline{\cal M}_g$, there is a universal  bundle
${\cal H}_g$, the Hodge line bundle.
Smith identifies the sum of the number of singular 
fibers and the signature to be $<4c_1({\cal H}_g), \phi_*[S^2]>$.

\medskip
\noindent{\bf \S3. Symplectic Parshin-Arakelov inequality}
\medskip

In this section, let $\pi:M\longrightarrow \Sigma$
be a genus $g$ relatively minimal Lefschetz fibration over a Riemann surface $\Sigma$.
 We will prove Theorem 1.

When $g=0$, $M$ is a ruled surface which is excluded by our assumption. 
When $g=1$, from the 
classification alluded before, $c_1^2=0$ and hence the inequality holds.
So let us assume $g\geq 2$. 
We can then (and will)
 choose a symplectic Lefschetz fibration structure on $M$ 
by Theorem 2.3.
First, we need to establish the following important fact.

\noindent{\bf Lemma 3.1}. There exist
 compatible 
almost complex structures on $M$
for which the fibers are pseudo-holomorphic submanifolds.

\noindent{\it Proof}. 
Near the singular point $x_i\in C$, the symplectic form
constructed by Gompf is K\"ahler with respect to
suitable local coordinates on which the projection has the form
$\pi(z_1, z_2)=z_1^2+z_2^2$. Fixing such an integrable complex structure $J_i$
in a closed neighborhood $U_i$ of each singular point, we see that
the intersection
$F_y\cap U_i$ is clearly holomorphic. 

Away from the singular point set $C$, the tangent bundle along 
the fibers $P$ is a symplectic sub-bundle. 
Its $\omega$ orthogonal dual $Q$ is also a symplectic subbundle. 
On the boundary of $(U_i)$, $P$ and $Q$ are both preserved by $J_i$.
It is well known that $J_i$ restricted to $P$ can be extended to 
a compatible complex structure on the complement of
$U_i$, and  the same is true for $Q$. Thus we obtain a 
compatible almost complex structure $J$ for which the fibers are pseudo-holomorphic. 

Let $F$ denote the class of fibers with complex orientation.

\noindent{\bf Proposition 3.2}. Suppose $M$ is not rational or ruled.
 Let $E$ be a class 
represented by an embedded sphere with square $-1$, which
has positive pairing with $\omega$,
 then $E\cdot F>0$.

\noindent{\it Proof}. Take a compatible almost complex structure $J$ 
constructed in the lemma above, the fibers are $J-$holomorphic
curves. 
For any compatible almost complex structure $J$,
 $E$ can be represented by a $J-$holomorphic curve $S$. 
This is true for $b^+>1$, as shown in [T].
In the case $b^+=1$, this follows from [LL1] with the additional 
assumption that $K\cdot E=-1$, and 
we ([L]) have proved that the assumption is always satisfied unless
$M$ is rational or ruled.  

If $E\cdot F\leq 0$, by the positivity of intersection,  $S$ must be 
contained in some singular fiber $F_s$, with its irreducible components
  also being irreducible components of $F_s$. This
is possible only if $F_s$ is a reducible fiber
 and one of its irreducible components
is a rational curve with square $-1$,  since we know $F_s$ has only one node.
 However, this contradicts with the assumption that $M\longrightarrow
\Sigma$ is  relatively minimal, and 
the proof is finished.

 Proposition 3.2  imply the following result of
Stipsicz. 

\noindent{\bf Corollary 3.3} ([S1]).
Suppose $\Sigma$ has positive genus and $ M\longrightarrow \Sigma$ is a relatively minimal 
Lefschetz fibration, then $M$ is minimal.

\noindent {\it Proof}. Observe that the intersection  number of a surface $S$  
with any fiber is simply the
degree of the restriction of the projection $\pi:S\longrightarrow 
\Sigma$, and when $S$ is a sphere and the genus of $\Sigma$ is positive,
 the degree has to be zero. 
If $M$ is not rational or ruled, it follows from Proposition 3.2 and the
observation that $M$ is minimal. 

We will finish the proof by showing that $M$ can not be rational or ruled.  
Suppose $M$ is rational or ruled. Choose a compatible almost complex structure
$J$ as constructed in Lemma 3.1. By [LL1], there exists an irreducible
$J-$holomorphic sphere $C$ representing a class $G$ with non-negative square. 
Since the fibers are $J-$holomorphic, $G\cdot F$ is non-negative.
$G\cdot F=0$ implies that $G$ is a irreducible component of
a singular fiber. But this is impossible because any irreducible
component of a singular fiber has square $-1$. 
$G\cdot F$ can not be positive either by the observation above. 
Thus the proof of Corollary 3.3 is finished.

We now prove Theorem 1 for  the cases $h\geq 1$ and $h=0$ in Theorems
3.4 and 3.5 separately. 

\noindent{\bf Theorem 3.4}. Let $M\longrightarrow \Sigma$ be a 
 relatively minimal Lefschetz fibration with fiber
$F$. If $g(\Sigma)\geq 1$ and $g(F)\geq 1$, then
$c_1^2(M)\geq 2(g(F)-1)(g(\Sigma)-1).$

\noindent{\it Proof}.   
Let $K$ denote the canonical class.
By Theorem 0.2 (1) in [T], $K$ is repesented by a smoothly embedded 
symplectic submanifold $C$. 
Furthermore, if $C_1, \cdots, C_k$ are the connected components of $C$, 
then for each $i$, $C_i^2\geq -1$.
If any $C_i$ is a sphere, since $M$ is minimal by Corollary 3.3, 
it must have non-negative self-intersection.
 This would imply that $M$ is rational or ruled according to
 a Theorem of McDuff ([Mc]), which is excluded by the claim in
the second paragraph of Corollary 3.3.  

%Then $b^+(M)$ is equal to one. But it was shown  that
%$b^+(M)$ is greater than three in [S1], hence none of the
%$C_i$ has genus zero. 

If we project $C_i$ to $\Sigma$, the degree of the 
projection is $d_i=C_i\cdot F$.
Since each $C_i$ has nonzero genus and $\Sigma$
has genus at least one, by a theorem of Kneser (see [Mi]), 
$$g(C_i)-1\geq d_i (g(\Sigma)-1).$$

Since the fibers are symplectic, we have 
the adjunction equality
$2(g(F)-1)=F\cdot F +K\cdot F=K\cdot F$.
Similarly, $2(g(C_i)-1)=C_i\cdot C_i+K\cdot C_i$.
Since $\sum_i (g(C_i)-1)=K^2$ and $\sum_j d_i=K\cdot F$,
we have 
$$K^2\geq 2(g(F)-1)(g(\Sigma)-1).$$
The theorem is proved becuse $c_1(M)$ is just $-K$.

\noindent {\bf Theorem 3.5}. Suppose $M$ is not rational or ruled and
$M\longrightarrow S^2$ is a  relatively minimal genus $g$
Lefschetz fibration over $S^2$,
then $c_1^2(M)\geq 2-2g$, and it is sharp.

\noindent{\it Proof}. 
Take a compatible almost complex structure $J$ constructed 
in Lemma 3.1. Let $E_1, \cdots, E_d$ be the
exceptional classes with $J$ holomorphic representatives 
$S_1\cdots, S_d$. 
Since $M$ is not ruled, according to [Mc], $M$ is obtained by
blowing up a minimal symplectic manifold $N$ at $d$ points with exceptional curves representing $E_1, \cdots, E_d$. Denote the blow-down map by $p$.
Thus $$K_M=p^*K_N+E_1+\cdots+E_d,\eqno (3.1)$$
and $K^2_M=K^2_N-d$.  
Since $K_N^2\geq 0$, by [T] and [Liu], it suffices to show that
$d\leq 2g-2$.
To prove $d\leq 2g-2$, we just need to show that
$$(E_1+\cdots+E_d)\cdot F\leq 2g-2,\eqno (3.2)$$
since $E_i\cdot F=S_i\cdot F\geq 1$.
By the adjunction formula,
$K\cdot F=2g-2$.
Thus (3.2) is equivalent to 
the claim that $p^*K_N\cdot F\geq 0$. 
This clearly  is true if $K_N$ is a torsion class.

Suppose $K_N$ is not a torsion class. When $b^+(M)>1$,
by Theorem 0.2 in [T] and by the blow-up formula of Gromov-Taubes
invariants in [LL2], the  class $p^*K_N$ can be 
represented by a $J$ holomorphic curve with components 
$T_1,\cdots, T_l$. 
 Since $F$ is represented by an irreducible $J$ holomorphic curve
with square zero, the positivity of intersection gives $T_i\cdot F\geq 0$ and hence the claim.
When $b^+(M)=1$,  $GT_N(2K_N)$ is shown to be nontrivial (see [LL3]).
Again by the blow-up formula of Gromov-Taubes invariants, $GT_M(p^*2K_N)=
GT_N(2K_N)$.
Thus $p^*2K_N$ can also be represented by $J$ holomorphic curve (possibly disconnected)
and we have the claim by similar argument. 

This inequality is sharp for many stable holomorphic Lefschetz
fibrations on blowups of K3 surfaces. Consider a generic pencil in a very 
ample system with square $2h$.  The base locus consists of $2h$ points.
The generic members are embedded curves with genus $h+1$ and the only 
singularity of each singular member is a nodal point.  Blow up the base locus,
we obtain a genus $h+1$ holomorphic stable Lefschetz fibration over
$S^2$ on $K3\#2h\overline{CP}^2$.  Clearly, $c_1^2$ and $2-2(1+h)$ are both 
equal to $-2h$.  Theorem 3.5 is proved.

It is proved in [S1] by a 
self fiber sum argument that $c_1^2\geq 4-4g$ for any 
Lefschetz fibrations. Examples of $M$ supporting Lefschetz fibrations
over $S^2$ with $c_1^2=4-4g$, as constructed in [GS],
are necessarily rational or ruled by Theorem 3.5.

Theorem 3.4 and 3.5 complete Theorem 1.

 Theorem 1 provides the symplectic analogue of Iitaka's conjecture 
$C_{2,1}$ concerning the Kodaira dimensions of the total space,
the fiber and the base of a stable holomorphic Lefschetz fibration.
We first introduce the definition of the symplectic analogue of the
Kodaira dimension.

\noindent{\bf Definition 3.6}.  The Kodaira dimension
$k(M)$ of a minimal symplectic $2-$manifold or a $4-$manifold with 
symplectic form $\omega$ and symplectic canonical class $K$  
is defined in the following way,

\item{} $k(M)=-\infty$ if $K\cdot \omega<0$;

\item{} $k(M)=0$ if $K\cdot \omega=0$;

\item{} $k(M)=1$ if $K\cdot \omega> 0$ and $K^2=0$;

\item{} $k(M)=2$ if $K\cdot \omega>0$ and $K^2>0$.

%\noindent{\bf Theorem} (Taubes) If $M$ is a minimal symplectic manifold with
%$b^+>1$, then $K\cdot K\geq 0$ and $K\cdot \omega\geq 0$.

%\noindent{\bf Theorem} (Liu) If $M$ is a mimimal symplectic manifold with
%$b^+=1$ and is not  ruled, then $K\cdot \omega \geq 0$.

The Kodaira dimension of a non-minimal symplectic $4-$manifold is  
the Kodaira dimension of one of its minimal models.

The definition for non-minimal symplectic $4-$manifolds does not depend on the choice of the minimal model.  
This is because that only rational and ruled symplectic four manifolds have 
more than one minimal models (see [Mc] and [L]), which 
 all have Kodaira dimension $-\infty$.

The symplectic manifolds with Kodaira dimension $-\infty$ have been 
classified. They are just the rational or ruled surfaces by results of
Taubes and Liu ([T], [Liu]).  We speculate that four-manifolds with Kodaira 
dimension zero either have K\"ahler structure or are torus bundles over torus.

With the above definition, the following is immediate from Theorem 1.

\noindent{\bf Corollary 3.7}. The Kodaira dimension of a Lefschetz fibration
is subadditive, ie. if $M\longrightarrow \Sigma$ is 
a Lefschetz fibration with fiber $F$, then $k(M)\geq k(F) + k(\Sigma)$.

\medskip
\noindent{\bf \S4. The number of singular fibers}
\medskip
In this section, let $M\longrightarrow S^2$ be
a Lefschetz fibration over $S^2$. We assume $g\geq 2$, since the cases for
$g=0$ and $g=1$ are well understood.
We have described in \S2 the correspondence between
Lefschetz fibrations over $S^2$ and relators consisting
of positive Dehn twists in the mapping class groups.
Given a Lefschetz fibration,  the number of singular fibers $l$ is 
just the length of the corresponding relator, and the number of irreducible 
singular fibers $n$ and the number of reducible sigular fibers $s$ 
are the numbers of positive Dehn twists along nonseparating 
curves and separating curves in the relator respectively. 
Here we study the lower bounds of $n$, $l$ and $s$.

The story for $s$ is very simple -- there are Lefschetz fibrations
with no reducible singular fibers for each $g$. In this section, we will 
focus on the lower bound of $n$.  We will also provide an estimate of 
the lower bound of $l$.

For the lower bound of $n$, we need to establish
the lower bound of the signature.

\noindent{\bf Proposition 4.1}. There are no Lefschetz fibrations over $S^2$
 with 
$\sigma=-l$.

\noindent{\bf Lemma 4.2}. Let $M\longrightarrow S^2$ be a
genus $g$ Lefschetz fibration on a ruled surface over a
genus $h$ Riemann surface $W$. Then 
$g\geq 2h-1$.

\noindent{\it Proof}.  When $h$ is 0, the statement is obvious.  
So we assume $h>0$.  
Let $F$ be a fiber. 
The composition of blowing down $p$ and projection $q$ to 
the base of the $S^2$-bundle $N$ gives rise to a smooth map 
$q\circ p:F\longrightarrow W$. 
It is shown in [LL1] that the fiber of any $S^2$-bundle has 
pseudo-holomorphic respresentative for any compatible almost 
complex structure. Thus the map $q\circ p$ must have positive 
degree because of the posivitity of intersection.
This implies $g\geq h$.
If we assume that $h\leq g<2h-1$, then any orientation-preserving map from
$\Sigma_g$ to $\Sigma_h$ must be of degree $1$.

Let us first assume that $N$ is the trivial $S^2-$bundle.
Let $U$ be a fiber and $V$ the section class of the $S^2-$bundle $N$
such that $U^2=V^2=0$ and $U\cdot V=1$.
Let $E_1, \cdots, E_k$ be the exceptional classes.
We denote the class of fibers of the Lefschetz fibration also by $F$. 
Since $U$, $V$, $E_1, \cdots, E_k$ form a basis of $H_2(M;{\bf Z})$,
$$F=aU+bV+c_1E_1+\cdots+c_kE_k$$
for some integers $a, b, c_1, \cdots, c_k,$ where $ c_i\leq 0$.
Since the degree of $q\circ p:F\longrightarrow W$ is $b$, 
$b$ must be $1$.  
From $F\cdot F=0$, we get $2a-c_1^2-\cdots-c_k^2=0$.
Recall that the canonical bundle is given
by $K=(2h-2)U-2V+E_1+\cdots +E_k$.
From the adjunction formula, we find $2h-2-2a-c_1-\cdots -c_k=2g-2$.
Thus $(c_1^2+c_1)+\cdots +(c_k^2+c_k)=2h-2g$.
Under the assumption that $h\leq g$, this is possible only if $h=g$.
However, in the case $h=g$, if $F$ is a reducible fiber, 
each of its component has genus less than $h$; if $F$ is an irreducible fiber,
its normalization has genus $h-1$. 
In either case, $b$ is forced to be $0$ rather than $1$, which 
leads to contradiction.

Similar argument applies to the case when $N$ is the nontrivial $S^2-$bundle.
The lemma is proved.

We now prove Proposition 4.1.

\noindent {\it Proof}.
Let $M\longrightarrow S^2$ be a genus $g$ Lefschetz fibration 
such that $\sigma=-l$.  It is easy to see that $b^+=1$ and $b_1=2g$.  
Then $M$ is the blow-up of a $S^2$-bundle over a genus $g$ surface by 
Theorem A in [Liu].  But under the assumption that $g\geq 2$, this contradicts Lemma 4.2.  
The proposition is proved.

Together with part 3 of Lemma 2.4, we immediately have the following corollary.

\noindent{\bf Corollary 4.3}. Any Lefschetz fibration over $S^2$ has at least one irreducible singular fiber.

In [ABKP], the authors conjectured that the monodromy group is not
contained in the Torelli group. This conjecture was proved in
[Sm2]. Since the Torelli group is generated by Dehn twists
along separating curves, their conjecture is also a consequence
of Corollary 4.3.

 Let $\mu(M)$ be the lowest genus of 
Lefschetz fibrations over $S^2$ on blowups of $M$. 
The $\mu$ invariant is zero  for ${CP}^2$ and $S^2\times S^2$, one for
elliptic surfaces,  and three for the four-torus  (see [Sm1]). 
In fact, with a little more effort we can determine $\mu$ for ruled surfaces. 

\noindent{\bf Proposition 4.4}. Let $M$ be a ruled surface over a
genus $h$ Riemann surface. Then $\mu= 2h$.

\noindent{\it Proof}.  When $h=0$, $\mu=2h$ is obvious. So we assume
$h\geq 1$. Suppose there is a genus $g=2h-1$ Lefschetz fibration.
By Kneser's theorem, the degree of $q\circ p$ is at most
two.
 On the other hand, by the argument in Lemma 4.2, we can rule out
the case when the degree of $q\circ p$
is one.  Therefore it must be exactly two. 
Consider an irreducible singular fiber,
whose existence is due to Corollary 4.3.
Its normalization is a surface of
genus $2h-2$, and therefore does not admit a degree two map to
a genus $h$ surface. Thus we have shown $g$ has to be greater than 
$2h-1$.

To show that $\mu=2h$, we need to construct
 $g=2h$  Lefschetz fibrations. There are many constructions
of such fibrations generalizing the example of Matsumoto ([M2]).
We sketch one here. 
Take  the trivial fibration $S^2\times \Sigma_h$, and let 
$U$ and $V$ be the fiber class and the section class as above. 
Consider the divisor class $U+2V$, which has square
four and its smooth members have
genus $g=2h$.  Blowing up four times, we obtain a genus
$2h$ Lefschetz fibrations. 

%$g\geq 2h+1$ when $g$ is odd. The equality is achieved 
%with fiber class $2U+2V$ on blowups at eight points.

\noindent{\bf Corollary 4.5}.
On any genus $g$ Lefschetz fibration over $S^2$ of the blowup of an
$S^2-$bundle, there are at least $2g$ singular fibers
and $g$ irreducible singular fibers.

% if $g$ is even, 
%and there are at least $2g+2$ singular fibers and $g+1$
%irreducble singular fibers if $g$ is odd.

\noindent{\it Proof}. Let $M$ be a blowup of an $S^2-$bundle 
over a surface of genus $h$. Suppose $M$ admits a Lefschetz 
fibration of genus $g$ with $l$ singular fibers.
Since the Euler number of an $S^2-$bundle over a surface
of genus $h$ is $-2(2h-2)$, we find that $l+-2(2g-2)\geq -2(2h-2).$
Since $g\geq 2h$, we have $l\geq 2g$.
 And by Lemma 2.4, we have $n\geq  2g-2h\geq g$. 

% if $h$ is even and $g\geq 2h+1$ if $h$
%is odd and the statement for the number of singular fibers follows.

\noindent{\bf Corollary 4.6}. Let $M\longrightarrow S^2$ be a 
genus $g$ Lefschetz fibration. If
$M$ is not rational or ruled, $n\geq (6g+6)/5+s/5$.

\noindent{\it Proof}.  By Theorem 3.5, we have $c_1^2\geq 2-2g$.
Since $c_1^2=2e+3\sigma$ and $e=4(g-1)(-1)+l$,
$\sigma= (-2l-8(g-1)(-1)+c_1^2)/3.$
By Lemma 2.5 and Proposition 4.1, $\sigma\leq n-s-4$, 
thus we find $n\geq (6g+6)/5+s/5.$

Now Theorem 2 follows from Corollaries 4.5 and 4.6.

We want to remark that when $g$ is odd, we can in fact show 
 $n\geq g+1$ with a more detailed 
analysis on ruled surfaces. A stronger bound 
for Lefschetz fibrations on
ruled surfaces is recently obtained in [S2].

When $g$ is low, it should be possible to determine the exact lower bound
of $n$.
We believe that the exact bound is six when $g=2$, and it is twelve
when $g=3$.  Examples with those numbers of irreducible singular fibers
include  genus
two Lefschetz fibration on $S^2\times T^2\# 4\overline{ CP}^2$ in [M2] 
and genus three fibrations on some torus bundles in [Sm1].  
We are not yet able to prove the exact bound, but we will present the best estimate in the following proposition.

% has exactly six irreducible singular
%fibers and two reducible singular fibers. 
%It is well known that $n$ is even ([M2]), so if $n$ is less  six, $n$ must be 
%four and $\sigma=-s=-3$.

\noindent{\bf Proposition 4.7}. The number of irreducible 
singular fibers of a genus $g$ Lefschetz fibration over $S^2$
is no less than four. And it is no less than six if $g\geq 3$.

\noindent{\it Proof}. By Lemma 2.4 and Lemma 2.5, 
$-n-s+4\leq \sigma\leq n-s-4$. So there are at least four 
irreducible singular fibers. If $g\geq 3$, the statement 
follows similarly from Lemma 2.5 and the following lemma.

\noindent{\bf Lemma 4.8}.
If $\sigma(M)=-l+4$, then $g\leq 2$. And 

\noindent 1. if $g=1$, $M$ is the rational elliptic
surface $E(1)$;

\noindent 2. if $g=2$, $M$ has $b^+=1$ and $b_1=2$.

\noindent{\it Proof}.
Let us assume that $\sigma=-l+4$. By Lemma 2.5, 
$b^+=1$ and $b_1=2g-2$. The last statement of the lemma follows.
If $g\geq 3$,  $M$ is the blow up
of an $S^2-$bundle over a genus $g-1$ surface according to [Liu].
But when $g\geq 3$, $g\leq 2(g-1)-1=2g-3$, which  is ruled out
by Proposition 4.4. 

If $g=1$,  we know $M$ is diffeomorphic to $E(k)$ 
hence $\sigma=-8k$ and $l=12k$. Therefore $\sigma=4-l=4-12k$ 
implies that $k=1$, so $M$ is $E(1)$.

Now we state a lower bound of the number of singular fibers, which 
follows from Cor. 4.5 and 4.6.

\noindent{\bf Proposition 4.9}. The number of singular fibers
in a genus $g$ Lefschetz fibration over $S^2$ is at least $(6g+6)/5$.

The best estimate, due to Stipsicz ([S2]), is $l\geq 8/5g$. 
We believe that the optimal bound of $l$ is of the order $2g$. 
 In fact, Gompf conjectures 
that the Euler number of a symplectic manifold $M$ 
is non-negative if $M$ is not a blow-up of an $S^2-$bundle over 
a surface of genus at least two. If this conjecture holds, 
then it is easy to see that there are at least $4g-4$ singular 
fibers for any Lefschetz fibration on any manifold which is 
not an $S^2$ bundle. With Corollary 4.5, we will be able to 
conclude that there are at least $2g+2$ singular fibers if 
$g$ is odd and $2g$ singular fibers if $g$ is even.

Recall that at the end of \S2, we introduced the geometric 
approach viewing genus $g$ Lefschetz fibrations as isotopy classes 
of smooth maps from the two sphere to the Deligne-Mumford moduli 
space of curves $\overline{\cal M}_g$ which have transverse 
positive intersections with ${\cal C}$. On $\overline{\cal M}_g$, 
there is a universal  bundle ${\cal H}_g$, the Hodge line bundle. 
Following  from Smith's signature formula and Corollary 4.6, 
 we have a positive lower bound of $c_1({\cal H}_g)$ 
linearly in the genus.

\noindent{\bf Corollary 4.10}. Suppose a genus $g$ Lefschetz fibration
over $S^2$ corresponds to a smooth map  
$\phi:S^2\longrightarrow \overline{\cal M}_g$. $$<c_1({\cal H}_g), 
\phi_*[S^2]> \geq {1\over 12}l+{g-1\over 3} \geq {3g-2\over 6}.$$

%If $M$ is not
%rational or ruled,  
%$$<c_1({\cal H}_g), \phi_*[S^2]>\geq {1\over 12}l+{g-1\over 2}
%\geq {3g-2\over 5}.$$
%If $M$ is rational or ruled, then the inequality follows from 
%Proposition 5.5. 

Notice that for holomorphic Lefschetz fibration, the positivity
is obvious since ${\cal H}_g$ is an ample line bundle and 
$\phi$ is a holomorphic map. As remarked in [Sm1], this is 
not a purely homological statement. Since by Wolpert's 
([W]) computation of the homology  of $\overline{\cal M}_g$,
there are two dimensional homology classes which have 
positive intersections with all the components of ${\cal C}$ 
but not with $c_1({\cal H}_g)$.

In [ABKP], the authors ask whether the pairing is still
non-negative when 
the genus $h$ of the base surface is positive. 
From Theorem 1 we can similarly derive  
$$<c_1({\cal H}_g), \phi_*[\Sigma_h]> 
\geq -{1\over 2}(h-1)(g-1)+{1\over 12}l.$$
This provides an affirmative answer to their question when $h=1$.

\medskip

\bigskip
\noindent{\bf References}.
\medskip
\item{}[ABKP] J. Amoros, F. Bogomolov, L. Katzarkov, T. Pantev,
Symplectic Lefschetz fibrations with arbitrary fundamental groups.
\item{} [D] S. Donaldson, Lefschetz fibrations in symplectic
geometry. Doc. Math. J. DMV., Extra Volume ICMII (1998) 309-314.
\item{} [E] H. Endo, Meyer's signature and hyperelliptic fibrations, preprint.
\item{}[GS] R. Gompf and A. Stipsicz, A introduction to 4-manifolds and Kirby calculus,
book in preparation.
\item{}[K] A. Kas, On the handlebody decomposition associated to
a Lefschetz fibration, Pacific. J. Math. 89 (1980) 89-104.
\item{} [Ko] D. Kotschick, Signatures, monoples and mapping
class groups, Math. Reserach Letter 5 (1998) 227-235.
\item{} [L] T. J. Li, Smoothly embedded spheres in 
symplectic four manifolds, Proc. AMS. 127 (1999) 609-613.
%\item{} [L2] T. J. Li, Lefschetz fibrations and 
%symplectic four manifolds, preprint.
\item{} [Liu] A-K. Liu, Some new applications of the general 
wall crossing formula, Math. Research Letters, 3 (1996) 569-585.
%\item{} [LL1], T. J. Li and A. K. Liu, General wall crossing formula,
%Math. Research Letters, 2 (1995) 797-810. 
\item{} [LL1], T. J. Li and A. K. Liu, Symplectic structures on ruled
surfaces and a generalized adjunction inequality, Math. Research 
Letters 2(1995) 453-471.
\item{} [LL2], T. J. Li and A. K. Liu, On the equivalence between
SW and Gr in the case $b^+=1$, IMRN,  (1999) 335-345.
\item{} [LL3], T. J. Li and A. K. Liu, Symplectic submanifolds
in symplectic four manifolds, in preparation.
\item{} [M1] Y. Matsumoto, Diffeomorphism types of elliptic 
surfaces, Topology 25 (1986) 549-563.
\item{} [M2] Y. Matsumoto, Lefschetz fibrations of genus 
two$-$a topological approach, Proceedings of the
37th Taniguchi Symposium on Topology and Teichmuller Spaces, ed. S.
Kojima et al., World Scientific (1996) 123-148.
\item{} [Mc]. D. McDuff, The structure of rational
and ruled symplectic $4-$manifold, Jour. AMS. v.1. no.3. (1990), 679-710.
%\item {}[Mc2]  D. McDuff, The local behavior of
%holomorphic curves in almost complex $4-$manifolds, Jour. Diff. Geom.,
%34(1991), 143-164.
\item{} [Mi] J. W. Milnor, On the existence of a connection with
curvature zero, Comm. Math. Helv. 32 (1959) 215-223.
\item{} [O] B. Ozbagci, Signatures of Lefschetz fibrations,
preprint. 
\item{} [Sm1] I. Smith, Symplectic geomtry of Lefschetz fibrations,
Oxford thesis, 1998.
\item{} [Sm2] I. Smith, Lefschetz fibrations and the Hodge bundle, 
Geometry and Topology 3 (1999) 211-233.
\item{} [S1] A. Stipsicz, Chern numbers of certain Lefschetz fibrations,
%and Erratum to this paper, 
Proc. AMS. to appear.
\item{} [S2] A. Stipsicz, Singular fibers in Lefschetz fibrations
on manifolds with $b^+=1$, preprint.
\item{} [ST] B. Siebert and G. Tian, On hyperelliptic $C^\infty-$Lefschetz
fibrations of four manifolds, Commun. Contemp. Math. 1 (1999) no.2 255-280.
\item{}[T] C. Taubes, $SW\Rightarrow Gr$: From Seiberg-Witten
equations to pseudo-holomorphic curves, JAMS 9 (1996) 845-918.
\item{} [W] S. Wolpert, On the homology of the moduli space of
stable curves, Ann. of Math., 118 (1983) 491-523.

\medskip
Department of Mathematics, Princeton University, Princeton, NJ 08544

tli@math.princeton.edu

\end